\documentclass{amsart}
\usepackage{amscd}
\usepackage[T1]{fontenc}   

\newtheorem{theorem}{Theorem}[section]
\newtheorem{lemma}[theorem]{Lemma}
\newtheorem{proposition}[theorem]{Proposition}
\newtheorem{corollary}[theorem]{Corollary}
\theoremstyle{definition}
\newtheorem{definition}[theorem]{Definition}

\newtheorem{remark}[theorem]{Remark}
\numberwithin{equation}{section}

\begin{document}

\noindent                                             
\begin{picture}(150,36)                               
\put(5,20){\tiny{Submitted to}}                       
\put(5,7){\textbf{Topology Proceedings}}              
\put(0,0){\framebox(140,34){}}                        
\put(2,2){\framebox(136,30){}}                        
\end{picture}                                        
\vspace{0.5in}

\renewcommand{\bf}{\bfseries}
\renewcommand{\sc}{\scshape}
\vspace{0.5in}

\title[On quasitopological homotopy groups of Inverse Limit Spaces]%
{On quasitopological homotopy groups of Inverse Limit Spaces}

\author{T. Nasri}
\address{Department of Pure Mathematics\\ Ferdowsi University of Mashhad\\
P.O.Box 1159-91775, Mashhad, Iran.}
\email{ta\_na8@stu-mail.um.ac.ir}

\author{ B. Mashayekhy}
\address{Department of Pure Mathematics\\ Center of Excellence in Analysis on Algebraic Structures\\ Ferdowsi University of Mashhad\\
P.O.Box 1159-91775, Mashhad, Iran.}
\email{bmashf@um.ac.ir}

\author{ H. Mirebrahimi}
\address{Department of Pure Mathematics\\ Ferdowsi University of Mashhad\\
P.O.Box 1159-91775, Mashhad, Iran.}
\email{h\_mirebrahimi@um.ac.ir}


\begin{abstract}
The paper is devoted to study the behavior of quasitopological homotopy groups on inverse limit spaces. More precisely, we present some conditions under which the quasitopological homotopy group of an inverse limit space and especially a product space is a topological group. Finally, we give some conditions for countability of homotopy groups.
\end{abstract}

\maketitle

\section{\bf Introduction and Motivation}
Endowed with the quotient topology induced by the natural surjective map $q:\Omega^n(X,x)\rightarrow \pi_n(X,x)$, where $\Omega^n(X,x)$ is the $n$th loop space of $(X,x)$ with the compact-open topology, the familiar homotopy group $\pi_n(X,x)$ becomes a quasitopological group which is called the quasitopological $n$th homotopy group of the pointed space $(X,x)$, denoted by $\pi_n^{qtop}(X,x)$ (see \cite{B,Br,Bra,G1}).

It was claimed by Biss \cite{B} that $\pi_1^{qtop}(X,x)$ is a topological group. However, Calcut and McCarthy \cite{CM} and Fabel \cite{F2} showed that there is a gap in the proof of \cite[Proposition 3.1]{B}. The misstep in the proof is repeated by Ghane et al. \cite{G1} to prove that $\pi_n^{qtop}(X,x)$ is a topological group \cite[Theorem 2.1]{G1} (see also \cite{CM}).

Calcut and McCarthy \cite{CM} showed that $\pi_1^{qtop}(X,x)$ is a homogeneous space and more precisely, Brazas \cite{Br} mentioned that $\pi_1^{qtop}(X,x)$ is a quasitopological group in the sense of \cite{A}.

Calcut and McCarthy \cite{CM} proved that for a path connected and locally path connected space $X$,  $\pi_1^{qtop}(X)$ is a discrete topological group if and only if $X$ is semilocally 1-connected (see also \cite{Br}). Pakdaman et al. \cite{P1} show that for a locally $(n-1)$-connected space $X$, $\pi_n^{qtop}(X,x)$ is discrete if and only if $X$ is semilocally n-connected at $x$ (see also \cite{G1}). Also, they prove that the quasitopological fundamental group of every small loop space is an indiscrete topological group.

Fabel \cite{F2,F3} and Brazas \cite{Br} presented some spaces for which their quasitopological homotopy groups are not topological groups.
 Moreover, despite of Fabel's result \cite{F2} that says the quasitopological fundamental group of the Hawaiian earring is not a topological group, Ghane et al. \cite{G2} proved that the quasitopological $n$th homotopy group of an $n$-Hawaiian like space is a prodiscrete metrizable topological group, for all $n\geq 2$.
By an n-Hawaiian like space $X$ we mean the
natural inverse limit, $\displaystyle{\lim_{\leftarrow}(Y_i^{n},y_i^*)}$, where $(Y_i^{n},y_i^*)=\bigvee_{j\leq i}(X_j^{n},x_j^*)$
is the wedge of $X_j^{n}$'s in which $X_j^{n}$'s are $(n-1)$-connected, locally $(n-1)$-connected, semilocally $n$-connected, and
compact CW spaces. Hence, it seems interesting to find out when $\pi_n^{qtop}(X,x)$ is or is not a topological group.

The recent result inspires the authors to study the behavior of quasitopological homotopy groups on inverse limit spaces.
Let $I$ be a partially order set, and $\{(X_i,x_i),\varphi_{ij}\}_I$ be an inverse system of pointed topological spaces, where $\varphi_{ij}:(X_j,x_j)\rightarrow(X_i,x_i)$ is a pointed continuous map, for all $i\leq j$. Suppose $(X,x)=\displaystyle{\lim_{\leftarrow}(X_i,x_i)}$ is the inverse limit space of the above inverse system. Since $\pi_n^{qtop}(-)$ is a functor from the category of pointed topological spaces, $Top_*$, to the category of quasitopological groups, $qTopGrp$, there exists a natural continuous homomorphism $\beta_n^{qtop}:\pi_n^{qtop}(X,x)\rightarrow \displaystyle{\lim_{\leftarrow}\pi_{n}^{qtop}(X_i,x_i)}$.
 The main goal of the paper is to present some inverse limit spaces and especially some product spaces whose quasitopological $n$th homotopy groups are topological groups, for $n\geq 1$.

The paper is organized as follows.
The crucial misstep in proving that $\pi_n^{qtop}(X,x)$ is a topological group is that for the natural quotient map $q:\Omega^n(X,x)\rightarrow \pi_n(X,x)$, the map $q\times q$ can fail to be a quotient map \cite[Theorem 1]{F2}. In Section 2, we intend to obtain some conditions in which the product map $q\times q$ is a quotient map for some spaces. Using this fact, we can present a class of product spaces whose quasitopological $n$th homotopy groups are topological groups. One of the main results of Section 2 is as follows:

{\it Let $(X,x)=\prod_{i\in I}(X_i,x_i)$, where $X_i$'s are second countable spaces whose $n$th homotopy groups are countable and Hausdorff. Then the isomorphism $\pi_n^{qtop}(X,x)\cong\prod_{i\in I}\pi_n^{qtop}(X_i,x_i)$ holds in topological groups.}

  Note that Biss \cite[Proposition 5.2]{B} claimed the above result without any conditions, but his proof has the mentioned misstep.
 Another main result of Section 2 is as follows:

{\it Let $(X,x)=\displaystyle{\lim_{\leftarrow}(X_i,x_i)}$  be the inverse limit space of an inverse system
$\{(X_i,x_i),\varphi_{ij}\}_I$, where $I$ is countable. Suppose that $X_i$ is second countable, $\pi_n^{qtop}(X_i,x_i)$ is Hausdorff, for all $i\in I$, $\pi_n^{qtop}(X,x)$ is second countable and the map $\beta_n:\pi_n(X,x)\rightarrow \displaystyle{\lim_{\leftarrow}\pi_{n}(X_i,x_i)}$ is a group monomorphism. Then $\pi_n^{qtop}(X,x)$ is a topological group.}\\

One of the main conditions of the above result is assuming that $\beta_n:\pi_n(X,x)\rightarrow \displaystyle{\lim_{\leftarrow}\pi_{n}(X_i,x_i)}$ is a group monomorphism. Using \cite{EK,FG,FZ}, we point out some classes of spaces for which the map $\beta_n$ is a group monomorphism. Moreover, we recall conditions under which homotopy group functors commute with inverse
limits.

In Section 3, we intend to generalize some results of Conner and Lamereaux \cite{CL} on the countability of $\pi_1(X,x)$ to obtain some conditions on a space to have a countable homotopy group. For this, we show that some properties of a space
 can be transferred to its loop space.

\section{\bf quasitopological homotopy groups of Inverse Limits}
In this section, first we intend to obtain some conditions on a topological space $X$ under which $\pi_n^{qtop}(X,x)$ is a topological group. The crucial misstep in proving that $\pi_n^{qtop}(X,x)$ is a topological group is that for the natural quotient map $q:\Omega^n(X,x)\rightarrow \pi_n(X,x)$, the map $q\times q$ can fail to be a quotient map (see \cite[Theorem 1]{F2}). We are going to find some conditions under which the product of quotient maps, $q\times q$, is also a quotient map. Michael \cite{Mi} introduced a new class of quotient maps, called bi-quotient maps. A map $f:X\rightarrow Y$ is called a bi-quotient map if, whenever $y\in Y$ and $\mathcal{U}$ is a covering of $f^{-1}(y)$ by open subsets of $X$, then finitely many $f(U)$, $U\in \mathcal{U}$, cover some neighborhood of $y$ in $Y$ \cite[Definition 1.1]{Mi}. He showed that any product (finite or infinite) of bi-quotient maps is also a bi-quotient map \cite[Theorem 1.2]{Mi}. Thus it is sufficient to see under which conditions the quotient map $q$ is a bi-quotient map. We need the following interesting results of Michael \cite{Mi}.
\begin{lemma}\label{L1} (\cite[Corollary 3.5]{Mi}).
 Let $f:X\rightarrow Y$ be a quotient map, where $Y$ is Hausdorff and $X$ is second countable. Then $f$ is a bi-quotient map if and only if $Y$ is second countable.
\end{lemma}

A topological space $X$ is called a {\it $k$-space} if a set $A\subseteq X$ is closed whenever $A\cap K$ is closed in $K$, for every compact $K\subseteq X$. Note that locally compact spaces and first countable spaces are $k$-spaces \cite{Mi}.
\begin{lemma}\label{L3}(\cite[Theorem 1.5]{Mi}). If $f_i:X_i\rightarrow Y_i$, $i=1,2$, are quotient maps and $X_1$ and $Y_1\times Y_2$ are both Hausdorff $k$-spaces, then $f_1\times f_2$ is a quotient map.
\end{lemma}

We shall also need the following well-known result.
\begin{theorem}\label{L2}(\cite[Theorem 12.5.2]{Du}).
If $X$ is second countable and $Y$ is locally compact and second countable, then the function space $X^Y$ is second countable. In particular, if $X$ is second countable then $\Omega^n(X,x)$ is also second countable, for all $x\in X$ and $n\in \Bbb{N}$.
\end{theorem}

The following theorems are two of the main results of this section.
\begin{theorem}\label{T1}
Suppose that $X$ is second countable,
 $\pi_n^{qtop}(X,x)$ is Hausdorff and second countable. Then $\pi_n^{qtop}(X,x)$ is a topological group.
\end{theorem}
\begin{proof}
  Since $X$ is second countable, by Theorem \ref{L2} $\Omega^n(X,x)$ is second countable. Now, by Lemma \ref{L1} the quotient map $q$ is a bi-quotient map. Hence by \cite[Theorem 1.2]{Mi} the product map $q\times q$ is a bi-quotient map and so it is a quotient map which yields the result.
\end{proof}
 \begin{theorem}\label{met}
 Let $X$ be a metric space. If one of the following conditions holds, then
  $\pi_n^{qtop}(X,x)$ is a topological group.\\
  $(i)$ $\pi_n^{qtop}(X,x)$ is Hausdorff and first countable.\\
  $(ii)$ $\pi_n^{qtop}(X,x)$ is Hausdorff and locally compact.
 \end{theorem}
 \begin{proof}
 $(i)$ Since $X$ is metric, the loop space $\Omega^n(X,x)$ is metric by \cite[Theorem 1 of Chapter IV \S44.V]{K} and so it is Hausdorff and first countable. Now, by Lemma \ref{L3} the product map $q\times q$ is a quotient map which implies the result.\\
$(ii)$ similar to the part $(i)$.
 \end{proof}

Now, we want to present a class of product spaces whose quasitopological $n$th homotopy groups are topological groups. Let $(X,x)=\prod_{i\in I}(X_i,x_i)$ and consider the following commutative diagram:
\begin{equation}\label{dia}
\begin{CD}
\Omega^n(X,x)@>\phi>>\prod_{i\in I}\Omega^n(X_i,x_i)\\
@VV q V@V \prod_{i\in I} q_i VV\\
\pi_n^{qtop}(X,x)@>\beta_n^{qtop}>>\prod_{i\in I}\pi_{n}^{qtop}(X_i,x_i),
\end{CD}\end{equation}
where $\phi$ is the natural homeomorphism using the fact that loop spaces preserve products, $q_i:\Omega^n(X_i,x_i)\rightarrow \pi_n^{qtop}(X_i,x_i)$ is the natural quotient map, for all $i\in I$. Since the map $q$ is a quotient map and the map $\beta_n:\pi_n(X,x)\rightarrow \prod_{i\in I}\pi_{n}(X_i,x_i)$ is a group isomorphism, if we show that $\prod_{i\in I} q_i$ is also a quotient map, then $\beta_n^{qtop}:\pi_n^{qtop}(X,x)\rightarrow \prod_{i\in I}\pi_{n}^{qtop}(X_i,x_i)$ will be an isomorphism of quasitopological groups. We are going to study the case in which the quotient map $q_i$ is a bi-quotient map, for all $i\in I$. The following theorem is another main result of this section.

\begin{theorem}\label{T11}
Let $\{(X_i,x_i)\}_{i\in I}$ be a family of second countable spaces such that $\pi_n^{qtop}(X_i,x_i)$ is Hausdorff and second countable, for all $i\in I$. Then the map $\prod q_i:\prod_{i\in I}\Omega^n(X_i,x_i)\rightarrow\prod_{i\in I}\pi_{n}^{qtop}(X_i,x_i)$ is a quotient map.
\end{theorem}
\begin{proof}
Since $X_i$ is second countable, by Theorem \ref{L2}, $\Omega^n(X_i,x_i)$ is  also second countable, for all $i\in I$. Since $\pi_n^{qtop}(X_i,x_i)$ is Hausdorff and second countable, by Lemma \ref{L1} the quotient map $q_i:\Omega^n(X_i,x_i)\rightarrow\pi_n^{qtop}(X_i,x_i)$ is a bi-quotient map, for all $i\in I$. Hence by \cite[Theorem 1.2]{Mi} the product map $\prod q_i$ is a bi-quotient map and so it is a quotient map.
\end{proof}

 It is well-known that if $(X,x)=\prod_{i\in I}(X_i,x_i)$, then $\pi_n(X,x)\cong\prod_{i\in I}\pi_n(X_i,x_i)$. However, it is known that the functor $\pi_n^{qtop}$ does not preserve products in general. For example, consider the Hawaiian Earring space, HE. Fabel \cite{F2} proved that the product $q\times q:\Omega(HE)\times\Omega(HE)\longrightarrow\pi^{qtop}_1(HE)\times\pi^{qtop}_1(HE)$ is not a quotient map and hence the topology of $\pi_1^{qtop}(HE\times HE)$ is strictly finer than $\pi_1^{qtop}(HE)\times \pi_1^{qtop}(HE)$. Moreover, Fabel \cite{F3} showed that for each $n\geq 1$ there exists a compact, path connected, metric space $X$ such that multiplication is discontinuous in $\pi_n^{qtop}(X,x)$. Hence $\pi_n^{qtop}$ does not preserve the product $X\times X$.
\begin{corollary}\label{C1}
Let $(X,x)=\prod_{i\in I}(X_i,x_i)$, where $X_i$'s are second countable spaces whose quasitopological $n$th homotopy groups are second countable and Hausdorff. Then the isomorphism $\pi_n^{qtop}(X,x)\cong\prod_{i\in I}\pi_n^{qtop}(X_i,x_i)$ holds in topological groups.
\end{corollary}
\begin{proof}
The result holds from Theorems \ref{T1} and \ref{T11} and using diagram \eqref{dia}.
\end{proof}
\begin{remark}\label{T111}
Let $\{(X_i,x_i)\}_{i\in I}$ be a family of metric spaces such that $\pi_n^{qtop}(X_i,x_i)$ is Hausdorff and first countable or Hausdorff and locally compact, for all $i\in I$. Then by Theorem \ref{met} and a similar proof of Theorem \ref {T11} and Corollary \ref{C1}, we have the following results, respectively.\\
$(i)$ The map $\prod q_i:\prod_{i\in I}\Omega^n(X_i,x_i)\rightarrow\prod_{i\in I}\pi_{n}^{qtop}(X_i,x_i)$ is a quotient map.\\
$(ii)$ The isomorphism $\pi_n^{qtop}(X,x)\cong\prod_{i\in I}\pi_n^{qtop}(X_i,x_i)$ holds in topological groups.
\end{remark}
\begin{proposition}\label{C2}
Let $(X,x)=\prod_{i\in I}(X_i,x_i)$, where $X_i$'s are locally $(n-1)$-connected and semilocally n-connected at $x_i$. Then the isomorphism $\pi_n^{qtop}(X,x)\cong\prod_{i\in I}\pi_n^{qtop}(X_i,x_i)$ holds in topological groups.
\end{proposition}
\begin{proof}
Since $X_i$ is locally $(n-1)$-connected and semilocally n-connected at $x_i$, by \cite[Theorem 6.7]{P1} $\pi_n^{qtop}(X_i,x_i)$ is discrete. Hence $\prod_{i\in I} q_i:\prod_{i\in I} \Omega^n(X_i,x_i)\rightarrow \prod_{i\in I} \pi_n(X_i,x_i)$  is open and the result follows from diagram (\ref{dia}).
\end{proof}

Note that the above proposition holds if we replace the conditions locally $(n-1)$-connected and semilocally n-connected by locally $n$-connected (see \cite[Theorem 6.6]{P1}).
The examples mentioned after Theorem \ref{T11} show that the assumptions second countability of $\pi_n^{qtop}(X_i,x_i)$'s in Corollary  \ref{C1} and semilocally n-connectedness in Proposition \ref{C2} are essential.
\begin{theorem}\label{dis}
If $\pi_n^{qtop}(X,x)$ is discrete, then $X$ is semilocally n-connected at $x$.
\end{theorem}
\begin{proof}
The result holds by a similar argument of \cite[Theorem 5.1]{B}, \cite[Lemma 3.1]{CM} and \cite[Theorem 3.2]{G1}.
\end{proof}

The following result seems interesting.
\begin{corollary}\label{sim}
Let $X=\prod_{i\in I}X_i$, where the $X_i$'s are locally $n$-connected. Then the following statements hold.\\
$(i)$ If $X$ is locally $n$-connected, then all but finitely many of the $X_i$'s are $n$-connected.\\
$(ii)$ The space $X$ is semilocally n-connected if all but finitely many of the $X_i$'s are $n$-connected.
\end{corollary}
\begin{proof}
$(i)$ Since $X$ is locally $n$-connected, $\pi_n^{qtop}(X)$ is discrete \cite[Theorem 6.6]{P1}. By Proposition \ref{C2}, $\prod_{i\in I}\pi_n^{qtop}(X_i)$ is discrete. Since the $X_i$'s are locally $n$-connected, $\pi_n^{qtop}(X_i)$ is discrete. But a product of infinitely many discrete spaces having more than one point is not discrete. Thus we conclude that all but finitely many of the $X_i$'s are $n$-connected.\\
$(ii)$ Applying Proposition \ref{C2}, we conclude that $\prod_{i\in I}\pi_n^{qtop}(X_i)\cong\pi_n^{qtop}(X)$ is discrete. Hence the result holds by Theorem \ref{dis}.
\end{proof}

As an immediate consequence of the above corollary, $\prod S^n$ is semilocally $k$-connected, for all $k<n$ and it is not locally n-connected.
\begin{corollary}\label{T2}
Let $(X,x)=\displaystyle{\lim_{\leftarrow}(X_i,x_i)}$  be the inverse limit space of an inverse system
$\{(X_i,x_i),\varphi_{ij}\}_I$, where $I$ is countable. Suppose that $X_i$ is second countable, $\pi_n^{qtop}(X_i,x_i)$ is Hausdorff, for all $i\in I$, $\pi_n^{qtop}(X,x)$ is second countable and the map $\beta_n:\pi_n(X,x)\rightarrow \displaystyle{\lim_{\leftarrow}\pi_{n}(X_i,x_i)}$ is a group monomorphism. Then
  $\pi_n^{qtop}(X,x)$ is a topological group.
\end{corollary}
\begin{proof}
Since $X_i$'s are second countable and $I$ is countable, the space $X$ is second countable by \cite[Theorem 8.6.2]{Du}. Since $\beta_n^{qtop}$ is an injective continuous map and $\pi_n^{qtop}(X_i,x_i)$'s are Hausdorff, $\pi_n^{qtop}(X,x)$ is Hausdorff. Hence Theorem \ref{T1} implies the result.
\end{proof}

Fabel \cite{F2} proved that the quasitopological fundamental group of the Hawaiian earring is not a topological group. Hence the hypothesis of $\pi_n^{qtop}(X,x)$ being second countable is essential in Corollary \ref{T2}.

One of the main conditions of the above result is assuming that $\beta_n:\pi_n(X,x)\rightarrow \displaystyle{\lim_{\leftarrow}\pi_{n}(X_i,x_i)}$ is a group monomorphism. In the following we point out some classes of spaces for which the map $\beta_1$ is a group monomorphism.
\begin{remark}
Let $(X,x)=\displaystyle{\lim_{\leftarrow} (X_i,x_i)}$ be compact and the $X_i$'s be compact polyhedra. Using \cite[Remark 1]{FZ}, the natural homomorphism $\beta_1:\pi_1(X,x)\rightarrow\displaystyle{\lim_{\leftarrow} \pi_{1}(X_i,x_i)}$ is a monomorphism if one of the following conditions holds:
\begin{enumerate}
\item $X$ is either a 1-dimensional compact Hausdorff space or $X$ is a 1-dimensional separable metric space (see \cite[Corollary 1.2]{EK}).
\item $X$ is any subset of a closed surface $M^{2}$ (see \cite[Theorem 5]{FZ}).
\item $X$ is a tree of $n$-manifolds (well-balanced if $n=2$) (see \cite[Theorem 3.1]{FG}).
\end{enumerate}
\end{remark}

\begin{remark}\label{R1}
Let $(X,x)=\displaystyle{\lim_{\leftarrow}(X_i,x_i)}$ be the inverse limit space of an inverse system
$\{(X_i,x_i),\varphi_{ij}\}_I$ and consider the following commutative diagram:
\begin{equation}\label{dia2}\begin{CD}
\Omega^n(X,x)@>\phi>>\displaystyle{\lim_{\leftarrow}\Omega^n(X_i,x_i)}\\
@VV q V@V QVV\\
\pi_n^{qtop}(X,x)@>\beta_n^{qtop}>>\displaystyle{\lim_{\leftarrow}\pi_{n}^{qtop}(X_i,x_i)},
\end{CD}\end{equation}
where $\phi$ is the natural homeomorphism and $Q$ is the restriction of the product of the quotient
maps $q_i$. Since the map $q$ is a quotient map, if we assume that the map $\beta_n:\pi_n(X,x)\rightarrow \displaystyle{\lim_{\leftarrow}\pi_{n}(X_i,x_i)}$ is a group isomorphism and show that $Q$ is also a quotient map, then $\beta_n^{qtop}:\pi_n^{qtop}(X,x)\rightarrow \displaystyle{\lim_{\leftarrow}\pi_{n}^{qtop}(X_i,x_i)}$ will be an isomorphism of quasitopological groups. Now, if the map $\prod q_i$ is a quotient map (see Theorem \ref{T11} and Remark \ref{T111} part (i)), the $X_i$'s are Hausdorff and $\displaystyle{\lim_{\leftarrow}\Omega^n(X_i,x_i)}$ is a saturated subspace of $\prod_{i\in I}(X_i,x_i)$ with respect to $\prod q_i$, then $Q$ is a quotient map \cite[Theorem 22.1]{Mu}.
\end{remark}

One of the main conditions of Remark \ref{R1} is assuming that $\beta_n:\pi_n(X,x)\rightarrow \displaystyle{\lim_{\leftarrow}\pi_{n}(X_i,x_i)}$ is a group isomorphism. In the following, we recall some conditions on the $X_i$'s in which the map $\beta_n$ is a group isomorphism.

The following interesting result of Cohen \cite{C} is essential.
\begin{theorem}\label{ex}
Let $(X,x)=\displaystyle{\lim_{\leftarrow}(X_i,x_i)}$ be the inverse limit space of an inverse system
$\{(X_i,x_i),\varphi_{ij}\}_I$, where the $X_i^,$s are Hausdorff and the consecutive maps, $\varphi_{ij}$, have the covering homotopy property for compactly generated spaces. Then the following short exact sequence exists, for all $n\geq 1$.
\begin{equation}
0\rightarrow \displaystyle{\lim_{\leftarrow}}^1 \pi_{n+1}(X_i,x_i)\rightarrow \pi_n(X,x)\overset {\beta_n}{\rightarrow} \displaystyle{\lim_{\leftarrow}} \pi_{n}(X_i,x_i)\rightarrow 0,
\end{equation}
where $\displaystyle{\lim_{\leftarrow}}^1$ is the first derived functor of the inverse limit functor.
\end{theorem}
\begin{remark}
To prove that $\beta_n$ is a group isomorphism, it is sufficient to see under which conditions $\displaystyle{\lim_{\leftarrow}}^1 \pi_{n+1}(X_i,x_i)=0$. With the assumptions of Theorem \ref{ex}, the natural homomorphism $\beta_n:\pi_n(X,x)\rightarrow\displaystyle{\lim_{\leftarrow} \pi_{n}(X_i,x_i)}$ is a group isomorphism,  for all $n\geq 1$ if $\{(X_i,x_i)\}$ is a movable inverse system. Indeed, movability of $X$ implies that $\displaystyle{\lim_{\leftarrow} \pi_{n}(X_i,x_i)}$ is movable by \cite[Remark 6.1.1]{MS} and
therefore it has Mittag- Leffler property \cite[Corollary 6.2.4]{MS}. Now by \cite[Theorem 6.2.10]{MS} we have
$\displaystyle{\lim_{\leftarrow}}^1\pi_n(X_i,x_i)=0$.
\end{remark}

\section{ \bf Countability of Homotopy Groups}

In this section, we intend to generalize some results of Conner and Lamereaux \cite{CL} on the countability of $\pi_1(X,x)$. For this, we show that some properties of topological spaces
can be transferred from $X$ to the loop space $\Omega^{n}(X,x)$, for some $x\in X$.
\begin{lemma}\label{semi}
Let $X$ be locally (n-1)-connected space. If $\Omega^{(n-1)}(X,x)$ is semilocally simply connected at the constant (n-1)-loop $e_x$, then $X$ is semilocally n-connected at $x$, for all $n\geq 2$.
\end{lemma}
\begin{proof}
Let $n=2$. Since $\Omega(X,x)$ is semilocally simply connected, there exists a basic element $V$ of the constant loop $e_x$ such that the homomorphism
  $\phi:\pi_1(V,e_x)\rightarrow\pi_1(\Omega(X,x),e_x)$, induced by inclusion, is trivial, where $V=\bigcap_{j=1}^{m}<K_{j},U_{j}>$. Put
 $U=\bigcap_{j=1}^{m}U_{j}$ and consider the following commutative diagram:
\[\begin{CD}
\pi_1(V,e_x)@>\phi>>\pi_1(\Omega(X,x),e_x)\\
@AA\theta^* A@AA\theta A\\
\pi_2(U,x)@>\psi>>\pi_2(X,x),
\end{CD}\]
where $\theta:\pi_2(X,x)\rightarrow \pi_1(\Omega(X,x),e_x)$ is given by $\theta([g])=[\bar{g}]$, where $\bar{g}(t)(s)=g(t,s)$ and $\theta^*$ is its restriction on $U$. Since $\theta:\pi_2(X,x)\rightarrow \pi_1(\Omega(X,x),e_x)$ is an isomorphism and the homomorphism $\phi$ is trivial, the homomorphism
 $\psi:\pi_2(U,x)\rightarrow\pi_2(X,x)$ is trivial. Hence $X$ is semilocally 2-connected at $x$. The result holds by induction on $n\geq 2$ similarly.
\end{proof}

Note that the converse of this fact has been shown by Wada \cite[Remark]{W}.
The first result on countability of homotopy groups is as follows.
\begin{theorem}\label{cont}
Let $X$ be an (n-1)-connected, locally (n-1)-connected, semilocally $n$-connected and separable metric space.
Then $\pi_n(X,x)$ is countable.
\end{theorem}
\begin{proof}
 It is easy to see that if $X$ is (n-1)-connected, then $\Omega^{(n-1)}(X,x)$ is path connected. Wada \cite[Corollary]{W} proved that if $X$ is locally (n-1)-connected, then $\Omega^{(n-1)}(X,x)$ is locally path connected and so it is locally connected. Also, he showed that the (n-1)-loop space of a semilocally n-connected space, is semilocally simply connected \cite[Remark]{W}. Since $X$ is a separable metric space, so is $\Omega^{(n-1)}(X,x)$ \cite[Theorem of Chapter II \S22.III]{K}. Therefore, $\Omega^{(n-1)}(X,x)$ satisfies the hypotheses of  \cite[Lemma 2.2]{CL} which implies that $\pi_1(\Omega^{(n-1)}(X,x),e_x)\cong\pi_n(X,x)$ is countable.
\end{proof}

A space $X$ is called $n$-homotopically Hausdorff at $x\in X$ if for
any essential n-loop $\alpha$ based at $x$, there is an open neighborhood $U$ of $x$ for
which $\alpha$ is not homotopic (rel $\dot{I^n}$) to any n-loop lying entirely in $U$.
$X$ is said to be n-homotopically Hausdorff if it is n-homotopically Hausdorff at
any $x\in X$ (see \cite{GH}).

Consider $\overline{\Omega^n(X,x)}$ as the space of homotopy classes rel $\dot{I^n}$ of n-loops at $x$ in $X$. If $p$ is an n-loop at $x$, and $U$ is an open neighborhood of $x$, then we define $O^n(p,U)$ to be the collection of homotopy classes of n-loops rel $\dot{I^n}$ containing n-loops of the form $p*\alpha$, where $\alpha$ is an n-loop in $U$ at $x$. It is routine to check that the collection $O^n(p,U)$ is a basis for $\overline{\Omega^n(X,x)}$. In the following, we show that $\overline{\Omega^n(X,x)}$ is Hausdorff if and only if $X$ is n-homotopically Hausdorff at $x$. ( see also \cite{CC} for the case $n=1$.)
\begin{lemma}\label{homo1}
$\overline{\Omega^n(X,x)}$ is Hausdorff if and only if $X$ is n-homotopically Hausdorff at $x$.
\end{lemma}
\begin{proof}
Let $f$ be an essential n-loop based at $x$. Since $\overline{\Omega^n(X,x)}$ is Hausdorff, there are neighborhoods $W=\bigcup_{i\in I} O^n(p_i,U_i)$ of $f$ and $W'=\bigcup_{i\in I} O^n(q_i,V_i)$ of $e_{x}$ such that $W\cap W'=\emptyset$. Put $V=\bigcup_{i\in I}V_i$, then $x\in V$. If $h$ is an n-loop in $V$ and $h\simeq f$ rel $\dot{I^n}$, then $[h]\in W\cap W'$ which is a contradiction. \\
Conversely, let $[f] , [g]\in \overline{\Omega^n(X,x)}$ such that $[f]\neq[g]$, then $g^{-1}*f\not\simeq e_{x}$. Since $X$ is n-homotopically Hausdorff, there is an open neighborhood $U$ of $x$ for
which $g^{-1}f$ is not homotopic (rel $\dot{I^n}$) to any n-loop lying entirely in $U$. Put $W=O^n(f, U)$ and $W'=O^n(g, U)$. We show that $W\cap W'=\emptyset$. Let $[h]\in W\cap W'$, then $[g*\alpha']=[h]=[f*\alpha]$, where $\alpha ,\alpha'$ are n-loops in $U$. Thus $[g^{-1}*f]=[\alpha'*\alpha^{-1}]$ which is a contradiction.
\end{proof}
\begin{lemma}\label{homo2}
$\overline{\Omega^n(X,x)}$ with the above topology is homeomorphic to\\ $\overline{\Omega(\Omega^{(n-1)}(X,x),e_x)}$, where $\Omega^{(n-1)}(X,x)$ is equipped with the compact-open topology, for all $n\geq 2$.
\end{lemma}
\begin{proof}
Consider the homeomorphism $\phi:\overline{\Omega^n(X,x)}\rightarrow\overline{\Omega(\Omega^{(n-1)}(X,x),e_x)}$ defined by $\phi([f])=[\overline{f}]$, where $\overline{f}(t)(s_1,s_2,\cdots,s_{n-1})=f(t,s_1,s_2,\cdots,s_{n-1})$.
\end{proof}
\begin{lemma}\label{homo} Let $n\geq 2$. Then
a space $X$ is n-homotopically Hausdorff at $x$ if and only if $\Omega^{(n-1)}(X,x)$ is homotopically Hausdorff at $e_x$, for any $x\in X$.
\end{lemma}
\begin{proof}
  By Lemma \ref{homo1}, $X$ is n-homotopically Hausdorff at $x$ if and only if $\overline{\Omega^n(X,x)}$ is Hausdorff. By Lemma \ref{homo2}, $\overline{\Omega^n(X,x)}$ is homeomorphic to
  $\overline{\Omega(\Omega^{(n-1)}(X,x),e_x)}$. Also, $\overline{\Omega(\Omega^{(n-1)}(X,x),e_x)}$ is Hausdorff if and only if $\Omega^{(n-1)}(X,x)$ is homotopically Hausdorff at $e_x$, for any
  $x\in X$.
\end{proof}

Now, the second result on countability of homotopy groups is as follows.
\begin{proposition}\label{uncont}
Suppose that $X$ is a second countable, locally (n-1)-connected and n-homotopically Hausdorff space at $x$ which is not semilocally n-connected at this point. Then $\pi_n(X,x)$ is uncountable.
\end{proposition}
\begin{proof}
This follows from \cite[Lemma 2.3]{CL}, Theorem \ref{L2} and Lemmas \ref{homo} and \ref{semi}.
\end{proof}

The following corollary is a consequence of Theorem \ref{cont} and Propositions \ref{uncont}.
\begin{corollary}
If $X$ is an (n-1)-connected, locally (n-1)-connected, separable metric space, then the following statements are equivalent.\\
$(i)$ $X$ is semilocally n-connected.\\
$(ii)$ $X$ is n-homotopically Hausdorff and $\pi_n(X)$ is countable.\\
\end{corollary}
\begin{definition}
Let $i : X \rightarrow Y$ be an embedding of one path connected
space into another. Then we say that $X$ is a $\pi_n$-retract of $Y$ if there
exists a homomorphism $r : \pi_n(Y ) \rightarrow \pi_n(X)$ such that the composition
$ri_* : \pi_n(X) \rightarrow \pi_n(X)$ is an isomorphism. In this case the homomorphism $r$ is called a $\pi_n$-retraction for X in Y.
Also, $X$ is called a $\pi_n$-neighborhood retract in $Y$ if $X$ is a $\pi_n$-retract of one of its open neighborhoods in $Y$.
\end{definition}
\begin{definition}
A separable metric space $X$ is called a $\pi_n$-absolute neighborhood
retract ($\pi_n$-$ANR$) if whenever $X$ is a subspace of
a separable metric space $Y$, then $X$ is a $\pi_n$-neighborhood retract in $Y$ .
\end{definition}
\begin{lemma}\label{semi2}
Let $Y$ be locally (n-1)-connected and semilocally n-connected and $X$ be a $\pi_n$-retract of $Y$. Then $X$ is semilocally n-connected.
\end{lemma}
\begin{proof}
Let $x\in X$. Since $Y$ is semilocally n-connected, there exists a neighborhood $V$ of $x$ in $Y$ such that $j_* : \pi_n(V,x) \rightarrow \pi_n(Y,x)$ is trivial. Since $X$ is a $\pi_n$-retract of $Y$, there
exists a homomorphism $r : \pi_n(Y ) \rightarrow \pi_n(X)$ such that the composition
$ri_* : \pi_n(X) \rightarrow \pi_n(X)$ is an isomorphism. We show that the induced mapping of inclusion $j'_* : \pi_n(V\cap X,x) \rightarrow \pi_n(X,x)$ is trivial. For this, consider the following commutative diagram:
\[\begin{CD}
 \pi_n(V\cap X,x)@>j'_*>>\pi_n(X,x)\\
@VV i'_* V@Vi_* VV\\
\pi_n(V,x)@>j_*>>\pi_n(Y,x),
\end{CD}\]
where the map $i'_*$ is induced by the inclusion. Now, the homomorphism $j_*$ is trivial and $i_*$ is injective, so $j'_*$ is trivial. Therefore $X$ is semilocally n-connected.
\end{proof}
\begin{corollary}
Let $X$ be a separable metric space. If $X$ is $\pi_n$-$ANR$, then it is semilocally n-connected.
\end{corollary}
\begin{proof}
Since X is a separable metric space, it follows from the proof of the Urysohn
metrization theorem \cite[Theorem 4.1]{Mu} that X can be embedded as a subspace of the
Hilbert cube $Q =\prod_{i=1}^{\infty}[0,1]$. Now, since $X$ is a $\pi_n$-$ANR$, we can choose an
open set, $U$, in the Hilbert cube such that $X$ is a $\pi_n$-retract in $U$. Since $U$ is semilocally n-connected, by Lemma \ref{semi2}, $X$ is also semilocally n-connected.
\end{proof}


\begin{proposition}
Let $X$ be an (n-1)-connected, locally (n-1)-connected, separable metric space in which $\pi_n(X,x)$ is free. Then $X$ is semilocally n-connected at $x$.
\end{proposition}
\begin{proof}
This follows from \cite[Theorem 2.6]{CL} and Lemma \ref{semi}.
\end{proof}

\bibliographystyle{plain}

\end{document}